\hsize 159.2mm
\vsize 246.2mm
\font\Bbb=msbm10
\font\bigrm=cmr17

\magnification=\magstep1
\def\C{\hbox{\Bbb C}}
\def\R{\hbox{\Bbb R}}

\footnote{}{ {\it Yang Xing

Department of Mathematics, University of Ume\aa, S-901 87 Ume\aa, Sweden

E-mail address:\enskip Yang.Xing@mathdept.umu.se
\bigskip

2000 Mathematics Subject Classification.} Primary  32W20, 32Q15 

 }

\centerline{\bigrm  The General Definition of the Complex Monge-Amp\`ere }
\smallskip
\centerline{\bigrm  Operator on Compact K\"ahler Manifolds }
\vskip .5in 
\centerline{\sl  Yang Xing      }
\vskip .8in
\item{}{{\bf Abstract.}\quad  We introduce a wide subclass ${\cal F}(X,\omega)$ of quasi-plurisubharmonic functions in a compact K\"ahler manifold, on which the complex Monge-Amp\`ere operator is well-defined and the convergence theorem is valid. We also prove that ${\cal F}(X,\omega)$ is a convex cone and includes all quasi-plurisubharmonic functions which are in the Cegrell class.
}
\bigskip\bigskip 

\noindent{\bf 1. Introduction }
\bigskip
Let $X$ be a compact connected K\"ahler manifold of dimension $n$, equipped with the fundamental form $\omega$ given in local coordinates by $\omega={i\over 2}\,\sum_{\alpha,\beta}g_{\alpha\bar\beta}dz^\alpha\wedge d\bar z^\beta,$
where $(g_{\alpha\bar\beta})$ is a positive definite Hermitian matrix and $d\omega=0$. 
The smooth volume form associated to this K\"ahler metric is the $n$th wedge product $\omega^n$.
Denote by $PSH(X,\omega)$ the set of upper semi-continuous functions $u:X\to\R\cup\{-\infty\}$ such that $u$ is integrable in $X$ with respect to the volume form $\omega^n$ and $\omega_u:=\omega+dd^cu\geq 0$ on $X$, where $d=\partial+\bar\partial$ and $d^c=i\,(\bar\partial-\partial)$.
These functions are called quasi-plurisubharmonic functions (quasi-psh for short) and play an important role in the study of positive closed currents in $X$, see Demailly's paper [D1]. A quasi-psh function is locally the difference of a plurisubharmonic function and a smooth function. Therefore, many properties of plurisubharmonic functions hold also for quasi-psh functions. Following Bedford and Taylor [BT2], the complex Monge-Amp\`ere operator $(\omega+dd^c)^n$ is locally and hence globally well defined for all bounded quasi-psh functions in $X$. Some important results of the complex Monge-Amp\`ere operator for bounded quasi-psh functions have been obtained by Kolodziej [KO1-2] and Blocki [BL1]. It is also known that the complex Monge-Amp\`ere operator does not work well for all unbounded quasi-psh functions. Otherwise, we shall lose some of the essential properties that the complex Monge-Amp\`ere operator should have, see Kiselman's paper [KI] or Bedford's survey [B]. In a bounded domain of $\C^n$ one usually needs certain assumptions on values of functions near the boundary of the domain to define complex Monge-Amp\`ere measures of unbounded plurisubahrmonic functions, see the Cegrell class [C1-2] where Cegrell introduced the largest subclass ${\cal E}(\Omega)$ of plurisuhharmonic functions in a bounded hyperconvex domain $\Omega$ for which the complex Monge-Amp\`ere operator is well-defined and the monotone convergence theorem is valid. However, such a technique does not seem to work for quasi-psh functions in a compact K\"ahler manifold because we lose boundary.
On the other hand, it was already observed by Bedford and Taylor [BT1] that for each quasi-psh function $u$ the complex Monge-Amp\`ere measure $\omega_u^n:=(\omega+dd^cu)^n$ is well defined on its non-polar subset $\{u>-\infty\}$. The complex Monge-Amp\`ere measures $\omega_u^n$ concentrating on $\{u>-\infty\}$ were studied by Guedj and Zeriahi [GZ]. In [X3] we obtained several convergence theorems for complex Monge-Amp\`ere measures without mass on pluripolar sets.
In this paper we introduce a quite large subclass ${\cal F}(X,\omega)$ of quasi-psh functions on which  images of the complex Monge-Amp\`ere operator are well-defined positive measures and may have positive masses on pluripolar sets. We prove that the set ${\cal F}(X,\omega)$ is a convex cone and includes all quasi-psh functions which are in the Cegrell class.
Our main result is the following convergence theorem of the complex Monge-Amp\`ere operator in ${\cal F}(X,\omega)$.
\bigskip
\noindent {\bf Theorem 5.}(Convergence Theorem) \it Let $0\leq p<\infty$. Suppose that $u_0\in {\cal F}(X,\omega)$ and that $g\in PSH(X,\omega)\cap L^\infty(X)$ is nonpositive. If $u_j,\,u\in {\cal F}(X,\omega)$ are such that $u_j\to u$ in $Cap_\omega$ on $X$ and $u_j\geq u_0$, then
$(-g)^p\,\omega_{u_j}^n\to (-g)^p\, \omega_u^n$ weakly in $X$.  
\rm
\bigskip
As a direct consequence we have
\bigskip
\noindent {\bf Corollary 5. } \it Let $0\leq p<\infty$ and $ 0\geq g\in PSH(X,\omega)\cap L^\infty(X)$. If $u_j,\,u\in {\cal F}(X,\omega)$ are such that $u_j\searrow u$ or $u_j\nearrow u$ in $X$, then
$(-g)^p\,\omega_{u_j}^n\to (-g)^p\,\omega_u^n$ weakly in $X$.  
\rm
\bigskip
For bounded quasi-psh functions, Corollary 5 is a slightly stronger version of the well-known monotone convergence theorem due to Bedford and Taylor [BT2].
\bigskip
\noindent{\it Acknowledgments}\quad I would like to thank Urban Cegrell for inspiring discussions on the subject. 
\bigskip
\bigskip 
\noindent{\bf 2.  The class ${\cal F}(X,\omega)$}
\bigskip
In this section we first introduce the subclass ${\cal F}(X,\omega)$ of quasi-psh functions, on which images of the complex Monge-Amp\`ere operator are finite positive measures in $X$. We obtain some characterizations of functions in ${\cal F}(X,\omega)$. Finally, we prove that ${\cal F}(X,\omega)$ is a star-shaped and convex set.

Recall that the Monge-Amp\`ere capacity $Cap_\omega$ associated to the K\"ahler form $\omega$ is defined by 
$$Cap_\omega(E)=\sup\bigl\{\int_E\omega_u^n;\,u\in PSH(X,\omega)\ {\rm and\ } -1\leq u\leq 0\bigr\},$$
for any Borel set $E$ in $X$. The capacity $Cap_\omega$ is introduced by Kolodziej [KO1] and is
comparable to the relative Monge-Amp\`ere capacity of Bedford and Taylor [BT2], and hence vanishes exactly on pluripolar sets of $X$.
Recall also that a sequence $\mu_j$  of positive Borel measures is said to be uniformly absolutely
continuous with respect to $Cap_\omega$ on $X$, or we write that $\mu_j\ll Cap_\omega$ on $X$ uniformly for all $j$, if for any 
$\varepsilon>0$ there exists $\delta>0$ such that 
$\mu_j(E)<\varepsilon$ for all $j$ and Borel sets $E\subset X$ with
$Cap_\omega(E)<\delta$.
Denote by $PSH^{-1}(X,\omega)$ the subset of functions $u$ in $PSH(X,\omega)$ with $\max\limits_Xu\leq -1.$
Given a function $u\in PSH^{-1}(X,\omega)$, we define the measure $(-u)\,\omega_u^{n-1}\wedge\omega$ in $X$ which is zero in $\{u=-\infty\}$ and 
$$\int\limits_E(-u)\,\omega_u^{n-1}\wedge\omega=\lim\limits_{j\to\infty}\int\limits_{E\cap\{u>-j\}}\bigl(-\max(u,-j)\bigr)\,\omega^{n-1}_{\max(u,-j)}\wedge\omega$$
for all $k\geq 1$ and $E\subset \{u>-k\}$. In a completely similar way, we define the measure $\omega_u^{n-1}\wedge\omega:=\chi_{\{u>-\infty\}}\,\omega_u^{n-1}\wedge\omega,$ where $\chi_{\{u>-\infty\}}$ is the characteristic function of the set $\{u>-\infty\}$. It is worth to point out that in general neither the measure $(-u)\,\omega_u^{n-1}\wedge\omega$ nor $\omega_u^{n-1}\wedge\omega$ is locally finite in $X$. However, we have the following result.
\bigskip
\noindent {\bf Proposition 1. } \it Let $u\in PSH^{-1}(X,\omega)$. Suppose that
$$-\max(u,-j)\,\omega^{n-1}_{\max(u,-j)}\wedge\omega \ll Cap_\omega\quad {\rm on}\  X\ {\rm uniformly\ for\ all}\ j=1,2,\dots.$$ 
Then the following statements hold.
\smallskip
\quad{\rm (1) $(-u)\,\omega_u^{n-1}\wedge\omega$ and $\omega_u^{n-1}\wedge\omega$ are finite positive measures in $X$;

\quad{\rm (2) $\max(u,-j)\,\omega^{n-1}_{\max(u,-j)}\to u\,\omega_u^{n-1}$ and $\omega^{n-1}_{\max(u,-j)}\to 
\omega_u^{n-1}$ as currents as $j\to\infty$;

\quad{\rm (3)  $(-u)\,\omega_u^{n-1}\wedge\omega\ll Cap_\omega$ on $X$.
\rm
\bigskip
\noindent{\it Proof.} Since 
$\int_X(-u)\,\omega_u^{n-1}\wedge\omega=\lim_{k\to\infty}\lim_{j\to\infty}\int_{u>-k}\bigl(-\max(u,-j)\bigr)\,\omega^{n-1}_{\max(u,-j)}\wedge\omega\leq \sup\limits_j\int_X\bigl(-\max(u,-j)\bigr)\,\omega^{n-1}_{\max(u,-j)}\wedge\omega<\infty,$ we obtain that $(-u)\,\omega_u^{n-1}\wedge\omega$ is a finite positive measure and so is $\omega_u^{n-1}\wedge\omega$.
Write $\max(u,-j)\,\omega^{n-1}_{\max(u,-j)}=\chi_{\{u\leq -j\}}\,\max(u,-j)\,\omega^{n-1}_{\max(u,-j)}$ $+\chi_{\{u> -j\}}\,\max(u,-j)\,\omega^{n-1}_{\max(u,-j)},$
where the first term on the right hand side tends to zero and the second one tends to $u\,\omega_u^{n-1}$ as $j\to\infty$. Similarly, we get that $\omega^{n-1}_{\max(u,-j)}\to 
\omega_u^{n-1}$  as $j\to\infty$.
Moreover, for any $E\subset X$ with $Cap_\omega(E)\not=0$ we can take an open set $G$ in $X$ such that $E \subset G$ and
$Cap_\omega(G)\leq 2\,Cap_\omega(E).$
Then  $\int_{E}(-u)\,\omega_u^{n-1}\wedge\omega\leq\int_G((-u)\,\omega_u^{n-1}\wedge\omega\leq\limsup\limits_{j\to\infty}
\int_G\bigl(-\max(u,-j)\bigr)\,\omega^{n-1}_{\max(u,-j)}\wedge\omega,$ which implies that $(-u)\,\omega_u^{n-1}\wedge\omega\ll Cap_\omega$ on $X$ and the proof of Proposition 1 is complete.
\bigskip
Let ${\cal F}(X,\omega)$ be the subset of functions in $PSH^{-1}(X,\omega)$ which satisfy the hypotheses of Proposition 1. The complex Monge-Amp\`ere measure $\omega_u^n$ of a function $u$ in ${\cal F}(X,\omega)$ is defined by the sum  
$$\omega_u^n:=\omega\wedge\omega_u^{n-1}+dd^c(u\,\omega_u^{n-1}),$$
where the currents $u\,\omega_u^{n-1}$ and $\omega_u^{n-1}$ are the limits of two sequences $\max(u,-j)\,\omega^{n-1}_{\max(u,-j)}$ and $\omega^{n-1}_{\max(u,-j)}$ respectively. Locally using the inequality $\bigl(\omega+dd^c(\phi+u)\bigr)^n\geq n\,\omega_u^{n-1}\wedge\omega$, where $\omega=dd^c\phi$, we can easily see that $(-u)\,\omega_u^{n-1}\wedge\omega\ll Cap_\omega$ in $X$ for any $u\in PSH^{-1}(X,\omega)\cap L^\infty(X)$, where $L^\infty(X)$ denotes the set of bounded functions in $X$. Hence for bounded quasi-psh functions, our definition of the complex Monge-Amp\`ere operator coincides with Bedford's and Taylor's definition given in [BT2].
Denote by $L^1(X,\mu)$ the set of integrable functions in $X$ with respect to the positive measure $\mu$. Now we give a characterization of functions in ${\cal F}(X,\omega)$.
\bigskip
\noindent {\bf Theorem 1. } \it Let $u\in PSH^{-1}(X,\omega)$. Then $u\in {\cal F}(X,\omega)$ if and only if $$u\in L^1(X,\,\omega_u^{n-1}\wedge\omega),$$ where $\omega_u^{n-1}:=\lim_{j\to\infty}\omega_{\max(u,-j)}^{n-1}$ as currents  and $\omega_{\max(u,-j)}^{n-1}\wedge\omega \ll Cap_\omega$ on $X$ uniformly for  $j=1,2,\dots$.
\rm
\bigskip
\noindent{\it Proof.} We prove first the " only if " part. Assume that $u\in {\cal F}(X,\omega)$. By Proposition 1 we have that $\omega_{\max(u,-j)}^{n-1}\wedge\omega\leq \bigl(-\max(u,-j)\bigr)\,\omega_{\max(u,-j)}^{n-1}\wedge\omega\ll Cap_\omega$ on $X$ uniformly for all $j$, and $\omega_{\max(u,-j)}^{n-1}\longrightarrow \omega_u^{n-1}$. Hence, by the lower semi-continuity of $-u$, we get that $\int_X\bigl(-\max(u,-t)\bigr)\,\omega_u^{n-1}\wedge\omega\leq \limsup_{j\to\infty}\int_X\bigl(-\max(u,-j)\bigr)\,\omega^{n-1}_{\max(u,-j)}\wedge\omega<\infty$ for all $t\geq 1$. Thus, we have $u\in L^1(X,\,\omega_u^{n-1}\wedge\omega)$. Now we prove the " if " part. Observe that for any $k> 1$, by Proposition 4.2 in [BT1] we get  $\chi_{\{u>-k\}}\,\omega_u^{n-1}\wedge\omega=\lim_{j\to\infty} \chi_{\{u>-k\}}\omega_{\max(u,-j)}^{n-1}\wedge\omega=\lim_{j\to\infty} \chi_{\{\max(u,-k)>-k\}}\omega_{\max(u,-j)}^{n-1}\wedge\omega=\lim_{j\to\infty} \chi_{\{\max(u,-k)>-k\}}$ $\omega_{\max(u,-j,-k)}^{n-1}\wedge\omega= \chi_{\{u>-k\}}\,\omega_{\max(u,-k)}^{n-1}\wedge\omega$. Hence, for any Borel set $E\subset X$ and $k>1$, we have that
$\int_E\omega_u^{n-1}\wedge\omega\leq \int_{u<-k+1}\omega_u^{n-1}\wedge\omega+\int_{E\cap\{u>-k\}}\omega_{\max(u,-k)}^{n-1}\wedge\omega\leq \limsup_{j\to\infty} \int_{u<-k+1}\omega_{\max(u,-j)}^{n-1}\wedge\omega+\int_E\omega_{\max(u,-k)}^{n-1}\wedge\omega,$
where we have used that the set $\{u<-k+1\}$ is open. Since $\omega_{\max(u,-j)}^{n-1}\wedge\omega \ll Cap_\omega$ on $X$ uniformly for  $j$, we have
 $\omega_u^{n-1}\wedge\omega\ll Cap_\omega$ on $X$. It then follows from $u\in L^1(X,\,\omega_u^{n-1}\wedge\omega)$ that $(-u)\,\omega_u^{n-1}\wedge\omega\ll Cap_\omega$ on $X$. For any $j\geq k_1>1$ we get 
$$\int\limits_{u\leq -k_1} \bigl(-\max(u,-j)\bigr)\,\omega_{\max(u,-j)}^{n-1}\wedge\omega\leq j\,\int\limits_{u\leq -j} \omega_{\max(u,-j)}^{n-1}\wedge\omega+\int\limits_{-j<u\leq -k_1} (-u)\,\omega_u^{n-1}\wedge\omega$$
$$= j\,\int\limits_X \omega^n-j\,\int\limits_{u>-j} \omega_{\max(u,-j)}^{n-1}\wedge\omega
+\int\limits_{-j<u\leq -k_1} (-u)\,\omega_u^{n-1}\wedge\omega$$
$$\leq  j\,\int\limits_X \omega^n-j\,\int\limits_{u>-j} \omega_u^{n-1}\wedge\omega
+\int\limits_{u\leq -k_1} (-u)\,\omega_u^{n-1}\wedge\omega\leq \int\limits_{\{u\leq-j\}\cup\{{u\leq -k_1}\}}(-u)\, \omega_u^{n-1}\wedge\omega.$$
Hence, for any Borel set $E_1\subset X$ and $j\geq k_1>1$, we have  
$\int_{E_1} \bigl(-\max(u,-j)\bigr)\,\omega_{\max(u,-j)}^{n-1}\wedge\omega
\leq \int_{\{u\leq-j\}\cup\{{u\leq -k_1}\}}(-u)\, \omega_u^{n-1}\wedge\omega 
+k_1\,\int_{E_1\cap\{u>-k_1\}} \omega_{\max(u,-j)}^{n-1}\wedge\omega:=A_{k_1,j}+B_{k_1,j}.$
Given $\varepsilon>0$ take $k_\varepsilon>1$ and $j_\varepsilon>1$ such that $A_{k_\varepsilon,j}\leq \varepsilon$ for all $j\geq j_\varepsilon$. Since $\omega_{\max(u,-j)}^{n-1}\wedge\omega \ll Cap_\omega$ on $X$ uniformly for  all $j$, there exists $\delta>0$ such that $(j_\varepsilon+k_\varepsilon)\,\int_{E_1}\omega_{\max(u,-j)}^{n-1}\wedge\omega\leq \varepsilon$ for all $j$ and $E_1\subset X$ with $Cap_\omega(E_1)\leq \delta.$ Therefore, we have proved that  $\int\limits_{E_1} \bigl(-\max(u,-j)\bigr)\,\omega_{\max(u,-j)}^{n-1}\wedge\omega\leq 2\,\varepsilon$ holds for all $j$ and $E_1\subset X$ with $Cap_\omega(E_1)\leq \delta.$ So $u\in {\cal F}(X,\omega)$ and the proof of Theorem 1 is complete.
\bigskip
Suppose that $\Omega$ is a hyperconvex subset in $\C^n$. Cegrell [C2] introduced the largest subclass ${\cal E}(\Omega)$ of plurisuhharmonic functions in $\Omega$, for which the complex Monge-Amp\`ere operator is well-defined and the monotone convergence theorem is valid. Our next theorem says that ${\cal F}(X,\omega)$ includes all quasi-psh functions which are  in the Cegrell class.  Recall that a negative plurisubharmonic function $u$ in $\Omega$ is said to belong to ${\cal E}(\Omega)$ if for each $z_0\in \Omega$ there exist a neighborhood $U_{z_0}$ of $z_0$ and a decreasing sequence $u_j$ of bounded plurisubharmonic functions in $\Omega$, vanishing on the boundary $\partial\Omega$, such that $u_j\searrow u$ on $U_{z_0}$ and $\sup\limits_j\int_\Omega(dd^cu_j)^n<\infty$. Blocki proved in [BL2] that it is a local property to belong to ${\cal E}(\Omega)$, that is, if $\Omega=\cup_j\Omega_j$ then $u\in {\cal E}(\Omega)$ if and only if $u\big|_{\Omega_j}\in {\cal E}(\Omega_j)$ for each $j$.
We call $u$ in $PSH^{-1}(X,\omega)$ for a Cegrell function in $X$ if there exists a finite covering $\{B_s\}_1^m$ of $X$ with hyperconvex subsets $B_s$ such that $\phi_s+u\in {\cal E}(B_s)$ for all $s$, where $\phi_s$ is a local K\"ahler potential defined in a  neighborhood of the closure of $B_s$, i.e. $\omega=dd^c\phi_s$ on $B_s=\{\phi_s<0\}$. Now we prove
\bigskip
\noindent {\bf Theorem 2. } \it If $u$ is a Cegrell function in $X$ then $u\in {\cal F}(X,\omega)$.
\rm
\bigskip
\noindent{\it Proof.} Take a new finite open covering $\{B_s^\prime\}_1^m$ of $X$ such that $B_s^\prime\subset\subset B_s$ for all $s$. By [C2] there exists a decreasing sequence $u_j^s$ of bounded plurisubharmonic functions in $B_s$, vanishing on  $\partial B_s$, such that $u_j^s\searrow \phi_s+u$ on $B_s^\prime$ and $\sup\limits_j\int_{B_s}(dd^cu_j^s)^n<\infty$. Since $Cap_\omega$ is
comparable to the relative Monge-Amp\`ere capacity of Bedford and Taylor, see [KO2][BT2], 
by Lemma 6 in [X2] we get that 
$-\max(u,-j)\,\omega^{n-1}_{\max(u,-j)}\wedge\omega \leq \bigl(-\phi_s-\max(u,-j)\bigr)\,\omega^{n-1}_{\max(u,-j)}\wedge\omega \ll Cap_\omega$ uniformly for all $j$ on 
each $B_s^\prime$ and hence on $X$. Therefore, $u\in {\cal F}(X,\omega)$ and the proof is complete.
\bigskip
Recall that a sequence $u_j$ of functions in $X$ is said to be convergent to a function $u$ in $ Cap_\omega$ on $X$ if for any $\delta>0$ we have
$$\lim\limits_{j\to\infty} Cap_\omega\bigl(\{z\in X;|u_j(z)-u(z)|>\delta\}\bigr)=0.$$ 
For a uniformly bounded sequence in $PSH(X,\omega)$, the convergence in capacity implies weak convergence of the complex Monge-Amp\`ere measures [X1]. 
Now we prove that the set ${\cal F}(X,\omega)$ is a convex cone. First, we need a lemma.
\bigskip
\noindent {\bf Lemma 1. } \it Let $u,\,v\in {\cal F}(X,\omega)$. Then
$$\int\limits_{u<v}(v-u)\,\omega_v^{n-1}\wedge\omega\leq \int\limits_{u<v}(v-u)\,\omega_u^{n-1}\wedge\omega.$$
If furthermore $u$ and $v$ are bounded, then for all integers $0\leq l\leq n-1$ we have 
$$\int\limits_{u<v}(v-u)\,\omega_v^l\wedge\omega_u^{n-1-l}\wedge\omega\leq \int\limits_{u<v}(v-u)\,\omega_u^{n-1}\wedge\omega.$$
\rm
\bigskip
\noindent{\it Proof.} We only prove the first inequality since the proof of the second one is similar. Assume first that $u$ and $v$ are bounded in $X$. By [D1] there exist a constant $A>1$ and two sequences $u_j,\, v_k\in PSH(X,A\omega)\cap C^\infty(X)$  such that $u_j\searrow u$ and $v_k\searrow v$ in $X$. Given $\varepsilon>0$. Assume first that $\{u_j<v_k\}\not=X$. Then $\max(v_k,u_j+\varepsilon)=u_j+\varepsilon$ near the boundary of the set $\{u_j<v_k\}$. Take a smooth subset $E_\varepsilon$ such that $\{u_j+\varepsilon<v_k\}\subset\subset E_\varepsilon\subset\subset\{u_j<v_k\}$, and write $T=\sum\limits_{l=0}^{n-2}\omega_u^l\wedge\omega_v^{n-2-l}\wedge\omega.$
Using Stokes theorem we get
$$\int\limits_{u_j<v_k}\bigl(\max(v_k,u_j+\varepsilon)-u_j-\varepsilon\bigr)\,\bigl((A\omega+dd^cu_j)-(A\omega+dd^c\max(v_k,u_j+\varepsilon))\bigr)\wedge T$$
$$=\int\limits_{E_\varepsilon}d\bigl(\max(v_k,u_j+\varepsilon)-u_j\bigr)\wedge d^c\bigl(\max(v_k,u_j+\varepsilon)-u_j\bigr)\wedge T\geq 0,$$
which holds even when $\{u_j<v_k\}=X$. Hence we obtain
$\int_{u_j<v_k}\bigl(\max(v_k,u_j+\varepsilon)-u_j\bigr)\,(A\omega+dd^cu_j)\wedge T\geq \int_{u_j<v_k}\bigl(\max(v_k,u_j+\varepsilon)-u_j-\varepsilon\bigr)\,(A\omega+dd^c\max(v_k,u_j+\varepsilon))\wedge T
\geq \int_{u_j<v_k}(v_k-u_j)\,(A\omega+dd^c\max(v_k,u_j+\varepsilon))\wedge T-\varepsilon\,A\,\int_X\omega^n.$
It turns out from the monotone convergence theorem in [BT2] that $(v_k-u_j)\,(A\omega+dd^c\max(v_k,u_j+\varepsilon))\wedge T\longrightarrow (v_k-u_j)\,(A\omega+dd^cv_k))\wedge T$ weakly in the open set $\{u_j<v_k\}$ as $\varepsilon\searrow 0$. Letting $\varepsilon\searrow 0$ and applying Lebesgue monotone convergence theorem we obtain the inequality
$\int_{u_j<v_k}(v_k-u_j)\,(A\omega+dd^cv_k)\wedge T\leq \int_{u_j<v_k}(v_k-u_j)\,(A\omega+dd^cu_j)\wedge T.$
Therefore, we have
$\int_{u_j<v}(v-u_j)\,(A\omega+dd^cv_k)\wedge T\leq \int_{u<v_k}(v_k-u)\,(A\omega+dd^cu_j)\wedge T.$
On the other hand, we have that $u_j,\,v_k$ are uniformly bounded, $u_j\to u$ in $Cap_\omega$ and $v_k\to v$ in $Cap_\omega$ on $X$. So for any $\delta>0$ the inequality 
$\int_{u<v}(v-u_j)\,(A\omega+dd^cv_k)\wedge T\leq \int_{u\leq v}(v_k-u)\,(A\omega+dd^cu_j)\wedge T+\delta$
holds for all $j,\, k$ large enough. Then by the quasicontinuity of quasi-psh functions, we can assume without loss of generality that $\{u<v\}$ is open and $\{u\leq v\}$ is closed. It turns out from the proof of Theorem 1 in [X1] that $(v-u_j)\,(A\omega+dd^cv_k)\wedge T\longrightarrow (v-u_j)\,(A\omega+dd^cu)\wedge T$ as $k\to\infty$ and $(v-u)\,(A\omega+dd^cu_j)\wedge T\longrightarrow (v-u)\,(A\omega+dd^cv)\wedge T$ as $j\to\infty$ weakly in $X$. Letting $k\to\infty$ and then $j\to\infty$, we obtain
$\int_{u<v}(v-u)\,(A\omega+dd^cv)\wedge T\leq \int_{u\leq v}(v-u)\,(A\omega+dd^cu)\wedge T+\delta.$
Applying $t\,v$ instead of $v$ for $A>t>1$ in the last inequality and then letting $t\searrow 1,\, \delta\searrow 0$ we get
$\int_{u<v}(v-u)\,(A\omega+dd^cv)\wedge T\leq \int_{u<v}(v-u)\,(A\omega+dd^cu)\wedge T,$
which yields that 
$\int_{u<v}(v-u)\,\omega_v^{n-1}\wedge\omega\leq \int_{u<v}(v-u)\,\omega_u^{n-1}\wedge\omega$
for all bounded quasi-psh functions $u$ and $v$. 

Now, for $u,\,v\in {\cal F}(X,\omega)$, we have
$\int_{\max(u,-j)<\max(v,-k)}\bigl(\max(v,-k)-\max(u,-j)\bigr)$ $\omega_{\max(v,-k)}^{n-1}\wedge\omega\leq \int_{\max(u,-j)<\max(v,-k)}\bigl(\max(v,-k)-\max(u,-j)\bigr)\,\omega_{\max(u,-j)}^{n-1}\wedge\omega.$
Letting $k\to\infty$, by the definition of $\omega_v^{n-1}\wedge\omega$ we get 
$\int_{\max(u,-j)<v}\bigl(v-\max(u,-j)\bigr)\,\omega_v^{n-1}\wedge\omega\leq \int_{\max(u,-j)<v}\bigl(v-\max(u,-j)\bigr)\,\omega_{\max(u,-j)}^{n-1}\wedge\omega,$
which by Fatou lemma implies that
$$\int\limits_{u<v}\bigl(v-u\bigr)\,\omega_v^{n-1}\wedge\omega\leq \liminf_{j\to\infty}\int\limits_{\max(u,-j)<v}\bigl(v-\max(u,-j)\bigr)\,\omega_{\max(u,-j)}^{n-1}\wedge\omega$$
$$\leq \liminf_{j\to\infty}\int\limits_{u<v}\bigl(\max(v,-j)-\max(u,-j)\bigr)\,\omega_{\max(u,-j)}^{n-1}\wedge\omega $$
$$\leq \limsup_{j\to\infty}\int\limits_{-s< u<v}\bigl(\max(v,-j)-\max(u,-j)\bigr)\,\omega_{\max(u,-j)}^{n-1}\wedge\omega $$
$$+\limsup_{j\to\infty}\int\limits_{\{u\leq-s\}\cap\{u<v\}}\bigl(\max(v,-j)-\max(u,-j)\bigr)\,\omega_{\max(u,-j)}^{n-1}\wedge\omega $$
$$=\int\limits_{-s< u<v}(v-u)\,\omega_u^{n-1}\wedge\omega+\limsup_{j\to\infty}\int\limits_{\{u\leq-s\}\cap\{u<v\}}\bigl(\max(v,-j)-\max(u,-j)\bigr)\,\omega_{\max(u,-j)}^{n-1}\wedge\omega $$
for all $s>1$. Since $\bigl(-\max(v,-j)\bigr)\,\omega_{\max(u,-j)}^{n-1}\wedge\omega< \bigl(-\max(u,-j)\bigr)\,\omega_{\max(u,-j)}^{n-1}\wedge\omega\ll Cap_\omega$ in the set $\{u<v\}$ uniformly for all $j$, letting $s\to\infty$ we get 
the required inequality and the proof of Lemma 1 is complete.
\bigskip
\noindent {\bf Theorem 3. } \it Let $u_0\in {\cal F}(X,\omega)$. If $u\in PSH^{-1}(X,\omega)$ satisfies $u\geq u_0$ in $X$ then
$u\in {\cal F}(X,\omega)$. Moreover, we have that $(-u)\,\omega_u^{n-1}\wedge\omega\ll Cap_\omega$ on $X$ uniformly for all  $u\in PSH^{-1}(X,\omega)$ with $u\geq u_0$ in $X$.
\rm
\bigskip
\noindent{\it Proof.} Given $k\geq 1$ and $j\geq 1$. Write $u_j=\max(u,-j)$. Then $u_j/3\in {\cal F}(X,\omega)$ and by Lemma 1 we have 
$\int_{u_j<-k}(-u_j)\omega_{u_j}^{n-1}\wedge\omega\leq 2  \int_{u_j<-k}(-k/2-u_j) \omega_{u_j}^{n-1}\wedge\omega\leq  3^{n-1}\,2\int_{u_j<-k/2}(-k/2-u_j) \omega_{{1\over 3} u_j}^{n-1}\wedge\omega\leq 3^n\int_{u_0<u_j/3-k/3}(u_j/3 -k/3-u_0)\,\omega_{{1\over 3} u_j}^{n-1}\wedge\omega\leq 3^n\int_{u_0<u_j/3-k/3}(u_j/3 -k/3-u_0)\,\omega_{u_0}^{n-1}\wedge\omega\leq 3^n\int_{u_0<-k/3}(-u_0)\,\omega_{u_0}^{n-1}\wedge\omega$. Thus, by $(-u_0)\,\omega_{u_0}^{n-1}\wedge\omega\ll Cap_\omega$ in $X$ we obtain that $(-u_j)\omega_{u_j}^{n-1}\wedge\omega\ll Cap_\omega$ in $X$ uniformly for all $j$, which yields that $u\in {\cal F}(X,\omega)$. Moreover, for all $k\geq 1$, $t\geq 1$ and $u\in PSH^{-1}(X,\omega)$ with $u\geq u_0,$ we have  
$\int_{\max(u,-t)<-k}(-u)\omega_u^{n-1}\wedge\omega \leq \limsup_{j\to\infty}\int_{\max(u,-t)<-k}(-u_j)$ $\omega_{u_j}^{n-1}\wedge\omega \leq \limsup_{j\to\infty}\int_{u_j<-k}(-u_j)\omega_{u_j}^{n-1}\wedge\omega \leq 3^n\int_{u_0<-k/3}(-u_0)\,\omega_{u_0}^{n-1}\wedge\omega$. Letting $t\to\infty$, we get $\int_{u<-k}(-u)\omega_u^{n-1}\wedge\omega\leq 3^n\int_{u_0<-k/3}(-u_0)\,\omega_{u_0}^{n-1}\wedge\omega$. Hence, together with $\chi_{\{u>-k-1\}}\,\omega_u^{n-1}\wedge\omega=\chi_{\{u>-k-1\}}\,\omega_{\max(u,-k-1)}^{n-1}\wedge\omega$, we obtain that $(-u)\,\omega_u^{n-1}\wedge\omega\ll Cap_\omega$ on $X$ uniformly for all  $u\geq u_0.$ The proof of Theorem 3 is complete.
\bigskip
As a direct consequence of Theorem 3 we have
\bigskip
\noindent {\bf Corollary 1. } \it Let $u\in {\cal F}(X,\omega)$. Then $\max(u,v)\in {\cal F}(X,\omega)$ and $t\,u\in {\cal F}(X,\omega)$ for all $v\in PSH^{-1}(X,\omega)$ and $0\leq t\leq 1$.
 \rm
\bigskip
Now we prove
\bigskip
\noindent {\bf Theorem 4. } \it The set ${\cal F}(X,\omega)$ is convex, that is, for any $u,\,v\in {\cal F}(X,\omega)$ and $0\leq t\leq 1$ we have that \quad$t\,u+(1-t)\,v\in {\cal F}(X,\omega)$.
 \rm
\bigskip
\noindent{\it Proof.} Given $u,\,v\in {\cal F}(X,\omega)$. Then $u/2+v/2\in PSH^{-1}(X,\omega)$. We only need to prove that $u/2+v/2\in {\cal F}(X,\omega)$. From Corollary 1 it turns out that $u/2\in {\cal F}(X,\omega)$ and $v/2\in {\cal F}(X,\omega)$. Then 
$\omega_{\max(u/2,-j)+\max(v/2,-j)}^{n-1}\wedge\omega=1/2^{n-1}\,\bigl(\omega_{\max(u,-2j)}+\omega_{\max(v,-2j)}\bigr)^{n-1}\wedge\omega\leq n!/2^{n-1}\,\sum\limits_{l=0}^{n-1}\omega_{\max(u,-2j)}^l\wedge \omega_{\max(v,-2j)}^{n-1-l}\wedge\omega.$ Write $u_{2j}=\max(u,-2j)$ and $v_{2j}=\max(v,-2j)$. For all $j\geq k\geq 1$ and $0\leq l\leq n-1$ we have  $$\int\limits_{u\leq -k}\omega_{u_{2j}}^l\wedge \omega_{v_{2j}}^{n-1-l}\wedge\omega= 1/k \,\int\limits_{u\leq -k}\bigl(-\max(u,-k)\bigr)\,\omega_{u_{2j}}^l\wedge \omega_{v_{2j}}^{n-1-l}\wedge\omega$$$$\leq 1/k\, \int\limits_X(-u_{2j})\,\omega_{u_{2j}}^l\wedge \omega_{v_{2j}}^{n-1-l}\wedge\omega\leq 1/k\,\int\limits_{u_{2j}\leq v_{2j}}  (-u_{2j})\,\omega_{u_{2j}}^l\wedge \omega_{v_{2j}}^{n-1-l}\wedge\omega$$$$+1/k \,\int\limits_{u_{2j}> v_{2j}}(-v_{2j})\,\omega_{u_{2j}}^l\wedge \omega_{v_{2j}}^{n-1-l}\wedge\omega.$$
From Lemma 1 it follows that $\int_{u_{2j}\leq v_{2j}}  (-u_{2j})\,\omega_{u_{2j}}^l\wedge \omega_{v_{2j}}^{n-1-l}\wedge\omega\leq 2\,\int_{u_{2j}\leq v_{2j}}  \bigl(v_{2j}/2-u_{2j}\bigr)$ $\omega_{u_{2j}}^l\wedge \omega_{v_{2j}}^{n-1-l}\wedge\omega\leq 2^{n-l}\,\int_{u_{2j}< v_{2j}/2}  \bigl(v_{2j}/2-u_{2j}\bigr)\,\omega_{u_{2j}}^l\wedge \omega_{v_{2j}/2}^{n-1-l}\wedge\omega\leq 2^{n-l}\,\int_{u_{2j}< v_{2j}/2}  \bigl(v_{2j}/2-u_{2j}\bigr)\,\omega_{u_{2j}}^{n-1}\wedge\omega\leq 2^{n-l}\,\sup\limits_j \int_X  (-u_{2j})\,\omega_{u_{2j}}^{n-1}\wedge\omega<\infty.$ Similarly, we have  
$\int_{u_{2j}> v_{2j}}(-v_{2j})\,\omega_{u_{2j}}^l\wedge \omega_{v_{2j}}^{n-1-l}\wedge\omega\leq 2^{l+1}\,\sup\limits_j \int_X  (-v_{2j})\,\omega_{v_{2j}}^{n-1}\wedge\omega<\infty.$ Hence we have proved that there exists a constant $A>0$ such that 
$\int_{\{u\leq -k\}\cup\{v\leq -k\}}\omega_{\max(u/2,-j)+\max(v/2,-j)}^{n-1}\wedge\omega\leq A/ k$ for all $j\geq k\geq 1$.
Thus, for $j\geq 2\,k\geq 1$ we have
$\int_{u/2+v/2\leq  -k}\omega_{\max(u/2+v/2,-j)}^{n-1}\wedge\omega=\int_X\omega^n-\int_{u/2+v/2>  -k}\omega_{\max(u/2+v/2,-j)}^{n-1}\wedge\omega=\int_X\omega^n-\int_{u/2+v/2>  -k}\omega_{\max(u/2,-j)+\max(v/2,-j)}^{n-1}\wedge\omega=\int_{u/2+v/2\leq  -k}\omega_{\max(u/2,-j)+\max(v/2,-j)}^{n-1}\wedge\omega\leq A/ k,$
which implies that $\omega_{\max(u/2+v/2,-j)}^{n-1}\wedge\omega\ll Cap_\omega$ on $X$ uniformly for all $j$ and hence $\omega_{u/2+v/2}^{n-1}\wedge\omega=\lim_{j\to\infty}\omega_{\max(u/2+v/2,-j)}^{n-1}\wedge\omega =\lim_{j\to\infty}\omega_{\max(u/2,-j)+\max(v/2,-j)}^{n-1}\wedge\omega.$
It then follows from the lower semi-continuity of $-u/2-v/2$ that $\int_X(-u/2-v/2)\,\omega_{u/2+v/2}^{n-1}\wedge\omega \leq \limsup_{j\to\infty}\int_X\bigl(-\max(u/4,-j/2)-\max(v/4,-j/2)\bigr)\,\omega_{\max(u/2,-j)+\max(v/2,-j)}^{n-1}\wedge\omega<\infty.$
By Theorem 1 we have obtained that $u/2+v/2\in {\cal F}(X,\omega)$, which concludes the proof of Theorem 4.
\bigskip
As consequences we have
\bigskip
\noindent {\bf Corollary 2. } \it Let $u_0,\,u_1,\dots,u_{n-1}\in {\cal F}(X,\omega)$. Then 
$$-u_0\,\omega_{u_1}\wedge\omega_{u_2}\wedge\dots\wedge\omega_{u_{n-1}}\wedge\omega\ll Cap_\omega\quad {\rm on }\quad X.$$
 \rm
\bigskip
\noindent{\it Proof.} Since $(u_0+u_1+\dots+u_{l-1})/l=(1/l)\,u_{l-1}+(1-1/l)\,(u_0+u_1+\dots+u_{l-2})/(l-1)$ for $l=2,3\dots,n$, using the induction principle and Theorem 4 we get that $f:=(u_0+u_1+\dots+u_{n-1})/n\in {\cal F}(X,\omega)$.
Hence we have that $-u_0\,\omega_{u_1}\wedge\omega_{u_2}\wedge\dots\wedge\omega_{u_{n-1}}\wedge\omega\leq -n^nf\,\omega_{u_1/n}\wedge\omega_{u_2/n}\wedge\dots\wedge\omega_{u_{n-1}/n}\wedge\omega\leq n^n\,(-f)\,\omega_f^{n-1}\wedge\omega \ll Cap_\omega$ on $X$, which concludes the proof of Corollary 2.
\bigskip
Using Corollary 2 and following the proof of Lemma 1, we get now a stronger version of Lemma 1.
\bigskip
\noindent {\bf Corollary 3. } \it Let $u,\,v\in {\cal F}(X,\omega)$ and $0\leq l\leq n-1$. Then 
$$\int\limits_{u<v}(v-u)\,\omega_v^l\wedge\omega_u^{n-1-l}\wedge\omega\leq \int\limits_{u<v}(v-u)\,\omega_u^{n-1}\wedge\omega.$$
\rm
\bigskip
\noindent {\bf Corollary 4. } \it Let $u_0\in {\cal F}(X,\omega)$. Then
$$-u_1\,\omega_{u_2}\wedge\omega_{u_3}\wedge\dots\wedge\omega_{u_n}\wedge\omega\ll Cap_\omega\quad {\rm on }\quad X$$
uniformly for all  $u_l\in PSH^{-1}(X,\omega)$ with $u_l\geq u_0$ and $l=1,2,\dots, n.$
\rm
\bigskip
\noindent{\it Proof.} Since $f:=(u_1+u_2+\dots+u_n)/n\geq u_0$ and $f\in {\cal F}(X,\omega)$, by Theorem 3 we get that $-u_1\,\omega_{u_2}\wedge\omega_{u_3}\wedge\dots\wedge\omega_{u_n}\wedge\omega\leq n^n\,(-f)\,\omega_f^{n-1}\wedge\omega \ll Cap_\omega$ on $X$
uniformly for all such functions $u_l$, which concludes the proof of Corollary 4.
\bigskip
\noindent{\it Remark.} Corollary 4 implies that a function $u\in PSH^{-1}(X,\omega)$ belongs to ${\cal F}(X,\omega)$ if and only if  $\bigl(-\max(u,-j)\bigr)\,\omega_{\max(u,-j)}^l\wedge\omega^{n-l}\ll Cap_\omega$ on $X$ uniformly for all $j\geq 1$ and $0\leq l\leq n-1$.

\bigskip
\bigskip 
\noindent{\bf 3.  A Convergence Theorem of the Complex Monge-Amp\`ere Operator}
\bigskip
In this section we prove a convergence theorem of the complex Monge-Amp\`ere operator in ${\cal F}(X,\omega)$. We divide its proof into several lemmas.

Given $u_1,\,u_2,\dots,u_{n-1}\in {\cal F}(X,\omega)$. By Corollary 2 the current $\omega_{u_1}\wedge\omega_{u_2}\wedge\dots\wedge\omega_{u_{n-1}}$ is well defined. Now for any $g\in PSH(X,\omega)\cap L^\infty(X)$, we define the wedge product 
$\omega_{u_1}\wedge\omega_{u_2}\wedge\dots\wedge\omega_{u_{n-1}}\wedge \omega_g$ 
in a natural way: $$\omega_{u_1}\wedge\omega_{u_2}\wedge\dots\wedge\omega_{u_{n-1}}\wedge \omega_g:=\omega\wedge\omega_{u_1}\wedge\omega_{u_2}\wedge\dots\wedge\omega_{u_{n-1}}+dd^c(g\,\omega_{u_1}\wedge\omega_{u_2}\wedge\dots\wedge\omega_{u_{n-1}}).$$ Then we have
\bigskip
\noindent {\bf Lemma 2. } \it Let $u_0,\,u_1,\dots,u_{n-1}\in {\cal F}(X,\omega)$ and $f,\,g\in PSH(X,\omega)\cap L^\infty(X)$. Then the following equalities hold.
$${\rm (a)}\qquad \int\limits_X(-g)\,dd^cf\wedge\omega_{u_1}\wedge\omega_{u_2}\wedge\dots\wedge\omega_{u_{n-1}}=\int\limits_X(-f)\,dd^cg\wedge\omega_{u_1}\wedge\omega_{u_2}\wedge\dots\wedge\omega_{u_{n-1}}.$$

$${\rm (b)}\qquad \int\limits_X(-g)\,dd^cu_0\wedge\omega_{u_1}\wedge\omega_{u_2}\wedge\dots\wedge\omega_{u_{n-1}}=\int\limits_X(-u_0)\,dd^cg\wedge\omega_{u_1}\wedge\omega_{u_2}\wedge\dots\wedge\omega_{u_{n-1}}.$$
\rm
\bigskip
\noindent{\it Proof.} It is no restriction to assume that $f,\,g\leq -2$ in $X$. Write $T=\omega_{u_1}\wedge\omega_{u_2}\wedge\dots\wedge\omega_{u_{n-1}}$. Take two sequences $f_j,\, g_k\in PSH^{-1}(X,A\omega)\cap C^\infty(X)$ for some $A\geq 1$ such that $f_j\searrow f$ and $g_k\searrow g$ in $X$, see [D1]. It follows from Dini's theorem and quasicontinuity of quasi-psh functions that $f_j\to f$ in $Cap_\omega$ on $X$.
So, using $T\wedge\omega\ll Cap_\omega$, 
we get $f_j\,T\to f\,T$ and hence $dd^cf_j\wedge T\to dd^cf\wedge T$ weakly in $X$. Similarly, $dd^cg_k\wedge T\to dd^cg\wedge T$ weakly in $X$. 
Thus we have
$\int_X(-f_j)\,dd^cg\wedge T=\lim_{k\to\infty}\int_X(-f_j)\,dd^cg_k\wedge T=\lim_{k\to\infty}\int_X(-g_k)\,dd^cf_j \wedge T=\lim_{k\to\infty}\int_X(-g_k)\,(A\omega+dd^cf_j) \wedge T-\lim_{k\to\infty}\int_X(-g_k)\,(A\omega) \wedge T=\int_X(-g)\,dd^cf_j \wedge T,$ where the last equality follows from the Lebesgue monotone convergence theorem.
Then, by lower semi-continuity of $-g$, we get   $\int_X(-f)\,dd^cg\wedge T=\lim_{j\to\infty}\int_X(-f_j)\,dd^cg\wedge T
=\lim_{j\to\infty}\int_X(-g)dd^cf_j \wedge T=\lim_{j\to\infty}\int_X(-g)(A\omega+dd^cf_j) \wedge T-\int_X(-g)(A\omega) \wedge T
\geq \int_X(-g)\,dd^cf \wedge T$. By symmetry we have abtained equality (a). Let $u_l=\max(u_0,-l)$. By (a) we have $\int_X(-g)\,dd^cu_l\wedge T=\int_X(-u_l)\,dd^cg\wedge T$. It follows from Corollary 2 that $u_0\,T$ is a well-defined current and $u_l\,T\to u_0\,T$ as currents in $X$. Hence we get $\int_X(-g)\,dd^cu_0\wedge T\leq lim_{l\to\infty}\int_X(-g)\,dd^cu_l\wedge T=lim_{l\to\infty}\int_X(-u_l)\,dd^cg\wedge T=\int_X(-u_0)\,dd^cg\wedge T.$
On the other hand, $\int_X(-u_0)\,dd^cg_k\wedge T=lim_{l\to\infty}\int_X(-u_l)\,dd^cg_k\wedge T=lim_{l\to\infty}\int_X(-g_k)\,dd^cu_l\wedge T=\int_X(-g_k)\,dd^cu_0\wedge T.$
Letting $k\to\infty$ we get $\int_X(-u_0)\,dd^cg\wedge T\leq \int_X(-g)\,dd^cu_0\wedge T$. Hence we have proved equality (b) and the proof of Lemma 2 is complete.
\bigskip
\noindent {\bf Lemma 3. } \it Let $u\in {\cal F}(X,\omega)$ and $g\in PSH(X,\omega)\cap L^\infty(X)$. Then the following  
statements hold.
\smallskip
\quad{\rm (a)} $\omega_{\max(u,-j)}^{n-1}\wedge \omega_g\ll Cap_\omega$ on $X$ uniformly for all  $j$;

\quad{\rm (b)} For each $f\in PSH(X,\omega)\cap L^\infty(X)$, we have that $f\,\omega_{\max(u,-j)}^{n-1}\wedge \omega_g\longrightarrow f\,\omega_u^{n-1}$

\qquad\quad  $\wedge \omega_g$ weakly in $X$ as $j\to\infty$;

\quad{\rm (c)} $(-u)\,\omega_u^{n-1}\wedge \omega_g\ll Cap_\omega$ on $X$.
\rm
\bigskip
\noindent{\it Proof.} It is no restriction to assume that $g\leq -2$ in $X$. Given $j\geq k\geq 1$. By Lemma 2 we have  $$\int\limits_{u\leq -k}\omega_{\max(u,-j)}^{n-1}\wedge \omega_g\leq 1/k\int\limits_X\bigl(-\max(u,-k)\bigr)\,\omega_{\max(u,-j)}^{n-1}\wedge \omega_g$$$$= 1/k\,\int\limits_X\bigl(-\max(u,-k)\bigr)\,\omega_{\max(u,-j)}^{n-1}\wedge \omega +1/k\,\int\limits_X(-g)\,\omega_{\max(u,-j)}^{n-1}\wedge dd^c\max(u,-k)$$
$$\leq 1/k\,\int\limits_X\bigl(-\max(u,-j)\bigr)\,\omega_{\max(u,-j)}^{n-1}\wedge \omega +1/k\,\int\limits_X(-g)\,\omega_{\max(u,-j)}^{n-1}\wedge \omega_{\max(u,-k)}$$
$$\leq 1/k\,\,\sup\limits_j\int\limits_X\bigl(-\max(u,-j)\bigr)\,\omega_{\max(u,-j)}^{n-1}\wedge \omega +1/k\,\sup\limits_X|g|\,\int_X\omega^n.$$
Given a Borel set $E\subset X$. By Proposition 4.2 in [BT1] for bounded quasi-psh functions, we get that $\int_E\omega_{\max(u,-j)}^{n-1}\wedge \omega_g\leq \int_{u\leq k}\omega_{\max(u,-j)}^{n-1}\wedge \omega_g+\int_E\omega_{\max(u,-k)}^{n-1}\wedge \omega_g$ for all $j\geq k\geq 1$, which implies (a).  

To prove (b), we prove first that $\omega_{\max(u,-j)}^{n-1}\wedge \omega_g\longrightarrow \omega_u^{n-1}\wedge \omega_g$ weakly in $X$ as $j\to\infty$. Given a smooth function $\psi$.  Multiplying a small positive constant if necessary, we can assume $\psi\in PSH(X,\omega)\cap C^\infty(X)$. Then we have
$\int_X\psi\,\omega_{\max(u,-j)}^{n-1}\wedge \omega_g-\int_X\psi\,\omega_u^{n-1}\wedge \omega_g=\int_X\psi\,\bigl(\omega_{\max(u,-j)}^{n-1}\wedge \omega-\omega_u^{n-1}\wedge \omega\bigr)+\int_Xg\,\bigl(\omega_{\max(u,-j)}^{n-1}-\omega_u^{n-1}\bigr)\wedge dd^c\psi,$
where by Proposition 1 the first term on the right hand side tends to zero as $j\to\infty$.
Take a sequence $ g_k\in PSH^{-1}(X,A\omega)\cap C^\infty(X)$ for some $A\geq 1$ such that $g_k\searrow g$ in $X$, see [D1]. Write the second term as
$$\int\limits_Xg_k\,\bigl(\omega_{\max(u,-j)}^{n-1}-\omega_u^{n-1}\bigr)\wedge dd^c\psi+\int\limits_X(g-g_k)\,\bigl(\omega_{\max(u,-j)}^{n-1}-\omega_u^{n-1}\bigr)\wedge dd^c\psi:=B_{k,j}+C_{k,j}.$$
By the smoothness of $\psi$ we have that $\bigl(\omega_{\max(u,-j)}^{n-1}+\omega_u^{n-1}\bigr)\wedge \omega_\psi\ll Cap_\omega$ on $X$ uniformly for all $j$. Since $g_k\to g$ in $Cap_\omega$ on $X$, we get that $C_{k,j}\to 0$ as $k\to\infty$ uniformly for all $j$. Then for each fixed $k$, $B_{k,j}\to 0$ as $j\to\infty$. Hence we have proved that 
$\omega_{\max(u,-j)}^{n-1}\wedge \omega_g\longrightarrow \omega_u^{n-1}\wedge \omega_g$ weakly in $X$ as $j\to\infty$.
Together with (a), we get $\omega_u^{n-1}\wedge \omega_g\ll Cap_\omega$ on $X$, see the proof of Proposition 1.
Now for $f\in PSH(X,\omega)\cap L^\infty(X)$, we take a sequence $ f_k\in PSH(X,A\omega)\cap C^\infty(X)$ for some $A\geq 1$ such that $f_k\searrow f$ in $X$. Write
$f\,\omega_{\max(u,-j)}^{n-1}\wedge \omega_g- f\,\omega_u^{n-1}\wedge \omega_g=(f-f_k)\,\bigl(\omega_{\max(u,-j)}^{n-1}\wedge \omega_g-\omega_u^{n-1}\wedge \omega_g\bigr)+f_k\,\bigl(\omega_{\max(u,-j)}^{n-1}\wedge \omega_g-\omega_u^{n-1}\wedge \omega_g\bigr),$ where for each fixed $k$ the second term on the right hand side tends to zero weakly as $j\to\infty$. Using (a) and $\omega_u^{n-1}\wedge \omega_g\ll Cap_\omega$, we get that the first term converges weakly to zero uniformly for all $j$ as $k\to\infty$. Thus we have obtained (b).

Finally, by the lower semi-continuity of $-u$, for any $k\geq 1$ we obtain  $\int_X\bigl(-\max(u,-k)\bigr)$ $\omega_u^{n-1}\wedge \omega_g\leq \limsup\limits_{j\to\infty}\int\limits_X\bigl(-\max(u,-k)\bigr)\,\omega_{\max(u,-j)}^{n-1}\wedge \omega_g
\leq \sup\limits_j\int\limits_X\bigl(-\max(u,-j)\bigr)\,\omega_{\max(u,-j)}^{n-1}$ $\wedge \omega +\sup\limits_X|g|\,\int_X\omega^n<\infty,$
which yields  $u\in L^1(X,\omega_u^{n-1}\wedge \omega_g)$. Thus we have that
$(-u)\,\omega_u^{n-1}\wedge \omega_g\ll \omega_u^{n-1}\wedge \omega_g\ll Cap_\omega$ on $X$. The proof of Lemma 3 is complete.
\bigskip
\noindent {\bf Lemma 4. } \it Let $u_0,\,u_1,\dots,u_{n-1}\in {\cal F}(X,\omega)$  and $g\in PSH(X,\omega)\cap L^\infty(X)$. Suppose that a sequence $u_j\in PSH^{-1}(X,\omega)$ decreases to $u_1$ in $X$. Then the following  
statements hold.
\smallskip
\quad{\rm (a)} $(-u_0)\,\omega_{u_1}\wedge\omega_{u_2}\wedge\dots\wedge\omega_{u_{n-1}}\wedge \omega_g\ll Cap_\omega$ on $X$;

\quad{\rm (b)} For each $f\in PSH(X,\omega)\cap L^\infty(X)$, we have that
$f\,\omega_{u_j}\wedge\omega_{u_2}\wedge\dots \wedge\omega_{u_{n-1}}\wedge \omega_g$ 

\qquad \quad  $\longrightarrow f\,\omega_{u_1}\wedge\omega_{u_2}\wedge\dots\wedge\omega_{u_{n-1}}\wedge \omega_g$ weakly in $X$ as $j\to\infty$;

\quad{\rm (c)} $\omega_{u_j}\wedge\omega_{u_2}\wedge\omega_{u_3}\wedge\dots\wedge\omega_{u_{n-1}}\wedge \omega_g\ll Cap_\omega$ on $X$ uniformly for all  $j$.
\rm
\bigskip
\noindent{\it Proof.} Since $(u_0+u_1+\dots+u_{n-1})/n\in {\cal F}(X,\omega)$, assertion (a) follows directly from (c) of Lemma 3. Now we prove (b). Given a smooth function $\psi$ in $X$. We assume without loss of generality that $0\leq f,\,\psi\in PSH(X,\omega)\cap L^\infty(X)$. Observe that $\varepsilon h^2\in PSH(X,\omega)$ if $h$ is a bounded positive quasi-psh function in $X$ and  the constant $\varepsilon$ satisfies $\max\limits_Xh\leq 1/(2\varepsilon)$.
Hence, applying the quality
${\psi\,f\over 2}=({\psi+f\over 2})^2-({\psi\over 2})^2-({f\over 2})^2$, we can assume that $h:=\psi\,f$ is a bounded quasi-psh function in $X$.  
By Lemma 2, for each $k\geq 1$ we get
$$\Bigl|\int\limits_X\psi\,f\,\omega_{u_j}\wedge\omega_{u_2}\wedge\dots \wedge\omega_{u_{n-1}}\wedge \omega_g-\int\limits_X\psi\,f\,\omega_{u_1}\wedge\omega_{u_2}\wedge\dots \wedge\omega_{u_{n-1}}\wedge \omega_g\Bigr|$$$$=\Bigl|\int\limits_X (u_j-u_1)\, dd^ch\wedge\omega_{u_2}\wedge\dots \wedge\omega_{u_{n-1}}\wedge \omega_g\Bigr|\leq \int\limits_X |u_j-u_1|\,(\omega_h+\omega)\wedge\omega_{u_2}\wedge\dots \wedge\omega_{u_{n-1}}\wedge \omega_g$$
$$\leq \int\limits_{u_1<-k} |u_1|\,(\omega_h+\omega)\wedge\omega_{u_2}\wedge\dots \wedge\omega_{u_{n-1}}\wedge \omega_g$$$$+\int\limits_X |\max(u_j,-k)-\max(u_1,-k)|\,(\omega_h+\omega)\wedge\omega_{u_2}\wedge\dots \wedge\omega_{u_{n-1}}\wedge \omega_g,$$
where by (a) the first term on the right hand side tends to zero as $k\to\infty$. For each fixed $k$, since $\max(u_j,-k)\to \max(u_1,-k)$ in $Cap_\omega$ on $X$ as $j\to\infty$, we get that the second term converges to zero as $j\to\infty$. Hence we have obtained (b).

By (a) and Theorem 3.2 in [BT1], assertion (c) follows from the property:  for any hyperconvex subset $\Omega\subset\subset X$ with $dd^c\phi=\omega$ and $\phi=0$ on $\partial\Omega$ and any $h\in PSH(\Omega)\cap L^\infty(\Omega)$, we have that $h\,\omega_{u_j}\wedge\omega_{u_2}\wedge\dots \wedge\omega_{u_{n-1}}\wedge \omega_g \longrightarrow h\,\omega_{u_1}\wedge\omega_{u_2}\wedge\dots\wedge\omega_{u_{n-1}}\wedge \omega_g$ weakly in $\Omega$ as $j\to\infty$. To prove this property, for each $\psi\in C^\infty_0(\Omega)$, we take a constant $\varepsilon>0$ such that $\varepsilon\,(h-\sup\limits_\Omega h-1)>\phi$ on $supp\, \psi$, and $\varepsilon\,(h-\sup\limits_\Omega h-1)<\phi$ near $\partial\Omega$. Set
$$ f=\cases{\max\bigl(\varepsilon\,(h-\sup\limits_\Omega h-1),\phi\bigr)-\phi,& in $\Omega$;\cr 0, & in $X\setminus \Omega$.\cr}. $$
Then $f\in PSH(X,\omega)\cap L^\infty(X)$ and $\psi\,h=\varepsilon^{-1}\psi\phi+\varepsilon^{-1}\psi f+\psi\, \sup\limits_\Omega h+\psi $. Hence, by the smoothness of $\phi$ and (b), we get that $h\,\omega_{u_j}\wedge\omega_{u_2}\wedge\dots \wedge\omega_{u_{n-1}}\wedge \omega_g \longrightarrow h\,\omega_{u_1}\wedge\omega_{u_2}\wedge\dots\wedge\omega_{u_{n-1}}\wedge \omega_g$ weakly in $\Omega$ as $j\to\infty$. Therefore, we have proved (c) and the proof of Lemma 4 is complete.
\bigskip
\noindent {\bf Lemma 5. } \it Let $u_0,\,u_1,\,u_2,\dots,u_{n-1}\in {\cal F}(X,\omega)$ and $g\in PSH(X,\omega)\cap L^\infty(X)$. Then for almost all constants $1\leq k<\infty$,
$$\int\limits_{u_1<-k}(-k-u_1)\,dd^cu_0\wedge\omega_{u_2}\wedge\dots \wedge\omega_{u_{n-1}}\wedge \omega_g\leq \int\limits_{u_1< -k}(-u_0)\,dd^cu_1\wedge\omega_{u_2}\wedge\dots \wedge\omega_{u_{n-1}}\wedge \omega_g.$$
\rm
\bigskip
\noindent{\it Proof.} Write $T=\omega_{u_2}\wedge\dots \wedge\omega_{u_{n-1}}\wedge \omega_g.$ Assume first that $0\geq u_0,\,u_1\in PSH(X,A\omega)\cap C^\infty(X)$ with $A\geq 1$. Given $\varepsilon>0$ and $k\geq 1$. Since $\max(u_1+\varepsilon,-k)=u_1+\varepsilon$ near $\partial\{u_1<-k\}$ if it is not empty, we have that $\int_{u_1<-k}(-k-u_1)\,dd^cu_0\wedge T=\lim_{\varepsilon\searrow 0}\int_{u_1<-k}\bigl(\max(u_1+\varepsilon,-k)-u_1-\varepsilon\bigr)\,dd^cu_0\wedge T=\lim_{\varepsilon\searrow 0}\int_{u_1<-k} u_0\,dd^c\bigl(\max(u_1+\varepsilon,-k)-u_1-\varepsilon\bigr)\wedge T=\int_{u_1<-k} (-u_0)\,dd^cu_1\wedge T+\lim_{\varepsilon\searrow 0}\int_{u_1<-k} u_0\,dd^c\max(u_1+\varepsilon,-k)\wedge T$. Since $\max(u_1+\varepsilon,-k)\, T\longrightarrow \max(u_1,-k)\, T$ weakly in $X$ as $\varepsilon\searrow 0$, we have $\bigl(A\omega+dd^c\max(u_1+\varepsilon,-k)\bigr)\wedge T\longrightarrow \bigl(A\omega+dd^c\max(u_1,-k)\bigr)\wedge T$ weakly as $\varepsilon\searrow 0$. From the upper semi-continuity of $u_0\leq 0$ in the open set $\{u_1<-k\}$, it turns out that $\lim_{\varepsilon\searrow 0}\int_{u_1<-k} u_0\,dd^c\max(u_1+\varepsilon,-k)\wedge T=\lim_{\varepsilon\searrow 0}\int_{u_1<-k} u_0\,\bigl[\bigl(A\omega+dd^c\max(u_1+\varepsilon,-k)\bigr)-A\omega\bigr]\wedge T
\leq \int_{u_1<-k} u_0\,dd^c\max(u_1,-k)\bigr)\wedge T
=0$. Hence we get $\int_{u_1<-k}(-k-u_1)\,dd^cu_0\wedge T\leq \int_{u_1<-k}(-u_0)\,dd^cu_1\wedge T$ for all $k\geq 1$ in the case of 
$0\geq u_0,\,u_1\in PSH(X,A\omega)\cap C^\infty(X)$.

Secondly, assume that $u_0,\,u_1\in {\cal F}(X,\omega)\cap L^\infty(X)$. By [D1] there exist negative functions $u_{0t},\,u_{1s}\in 
PSH(X,A\omega)\cap C^\infty(X)$ with some $A\geq 1$ such that $u_{0t}\searrow u_0$ and $u_{1s}\searrow u_1$ in $X$.
Since $\int_{u_1\leq -k}\bigl(\omega_{u_1}+\omega\bigr)\wedge T$ is an increasing function of $k$ and hence continuous almost everywhere with respect to the Lebesgue measure, we have that $\int_{u_1=-k}\bigl(\omega_{u_1}+\omega\bigr)\wedge T=0$ holds
for almost all $k$ in $[1,\infty).$ Given such a constant $k$. By Fatou lemma and the lower semi-continuity of $-u_{1s}$, we get that $\int_{u_1<-k}(-k-u_1)\,dd^cu_0\wedge T=\int_{u_1<-k}(-k-u_1)\,(A\omega+dd^cu_0)\wedge T-A\,\int_{u_1<-k}(-k-u_1)\,\omega\wedge T
\leq \liminf_{s\to\infty}\int_{u_{1s}<-k}(-k-u_{1s})\,(A\omega+dd^cu_0)\wedge T-A\,\int_{u_1<-k}(-k-u_1)\,\omega\wedge T
\leq \liminf_{s\to\infty}\limsup_{t\to\infty}\int_{u_{1s}<-k}(-k-u_{1s})\,(A\omega+dd^cu_{0t})\wedge T-\liminf_{s\to\infty}A\,\int_{u_1<-k}(-k-u_{1s})\,\omega\wedge T=\liminf_{s\to\infty}\limsup_{t\to\infty}\int_{u_{1s}<-k}(-k-u_{1s})\,dd^cu_{0t}\wedge T-A\,\liminf_{s\to\infty}\int_{u_{1s}\geq -k>u_1}(-k-u_{1s})\,\omega\wedge T.$ Given $\delta>0$, we have that
$\bigl|\int_{u_{1s}\geq -k>u_1}(-k-u_{1s})\,\omega\wedge T\bigr|\leq \delta\,
\int_X\omega\wedge T+\int_{u_{1s}-u_1\geq \delta}(-u_1)\,\omega\wedge T\longrightarrow \delta\,
\int_X\omega\wedge T$ as $s\to\infty$, since $u_{1s}\to u_1$ in $Cap_\omega$ and $(-u_1)\,\omega\wedge T\ll Cap_\omega$ on $X$. Hence we have 
$\int_{u_1<-k}(-k-u_1)\,dd^cu_0\wedge T\leq \liminf_{s\to\infty}\limsup_{t\to\infty}\int_{u_{1s}<-k}(-k-u_{1s})\,dd^cu_{0t}\wedge T\leq \liminf_{s\to\infty}\limsup_{t\to\infty}\int_{u_{1s}<-k}(-u_{0t})\,dd^cu_{1s}\wedge T =\liminf_{s\to\infty}\int_{u_{1s}<-k}(-u_0)\,dd^cu_{1s}\wedge T\leq \liminf_{s\to\infty}\int_{u_1\leq -k}(-u_0)\,(A\omega+dd^cu_{1s})\wedge T-A\,\liminf_{s\to\infty}\int_{u_{1s}<-k}(-u_0)\,\omega\wedge T=\liminf_{s\to\infty}\int_{u_1\leq -k}(-u_0)\,(A\omega+dd^cu_{1s})\wedge T-A\,\int_{u_1\leq -k}(-u_0)\,\omega\wedge T$. By Lemma 4 and quasicontinuity of quasi-psh functions, it is no restriction to assume that $\{u_1\leq -k\}$ is a closed set and hence the last limit inferior does not exceed $\int_{u_1\leq -k}(-u_0)\,(A\omega+dd^cu_1)\wedge T$. So we have obtained $\int_{u_1<-k}(-k-u_1)\,dd^cu_0\wedge T\leq \int_{u_1< -k}(-u_0)\,dd^cu_1\wedge T$
for all  $u_0,\,u_1\in {\cal F}(X,\omega)\cap L^\infty(X)$ and almost all $k$ in $[1,\infty).$

Finally, let $u_0,\,u_1\in {\cal F}(X,\omega)$. For almost all constants $k$ in $[1,\infty)$ we have that $\int_{u_1=-k}\bigl(\omega_{u_1}+\omega\bigr)\wedge T=0$ and $\int_{\max(u_1,-s)<-k}\bigl(-k-\max(u_1,-s)\bigr)\,dd^c\max(u_0,-t)\wedge T\leq \int_{\max(u_1,-s)< -k}\bigl(-\max(u_0,-t)\bigr)\,dd^c\max(u_1,-s)\wedge T$
for all integers $\,s,\,t\geq 1$. Letting $s\to\infty$ and applying the same proof as above, we have 
$\int_{u_1<-k}(-k_j-u_1)\,dd^c\max(u_0,-t)\wedge T\leq \int_{u_1< -k}\bigl(-\max(u_0,-t)\bigr)\,dd^cu_1\wedge T$
and then letting $t\to\infty$ we get the required inequality. The proof of Lemma 5 is complete.
\bigskip
\bigskip
\noindent {\bf Lemma 6. } \it Let $u_0\in {\cal F}(X,\omega)$ and $g\in PSH(X,\omega)\cap L^\infty(X)$. Then
$$\int\limits_{u<-k}(-u)\,\omega_u^{n-1}\wedge\omega_g\longrightarrow 0,\qquad{\rm as\quad} k\to\infty,$$
uniformly for all $u\in PSH^{-1}(X,\omega)$ with $u\geq u_0$ in $X$.
\rm
\bigskip
\noindent{\it Proof.} Given $u\in PSH^{-1}(X,\omega)$ with $u\geq u_0$. Take a sequence $1\leq k_1\leq k_2\leq \dots\leq k_j\to\infty$ such that Lemma 5 holds for the functions $u$ and $u_0$ when $k=k_j/2^i$, where $i=1,\dots,n-1$ and $j=1,2,\dots.$ Hence we have
$$\int\limits_{u<-k_j}(-u)\,\omega_u^{n-1}\wedge\omega_g\leq \int\limits_{u_0<-k_j}(-u_0)\,\omega_u^{n-1}\wedge\omega_g\leq 2\,\int\limits_{u_0<-k_j}(-k_j/2-u_0)\,\omega_u^{n-1}\wedge\omega_g$$$$\leq 2\,\int\limits_{u_0<-k_j/2}(-k_j/2-u_0)\,\omega\wedge\omega_u^{n-2}\wedge\omega_g+2\,\int\limits_{u_0<-k_j/2}(-k_j/2-u_0)\,dd^cu\wedge\omega_u^{n-2}\wedge\omega_g$$$$\leq 2\,\int\limits_{u_0<-k_j/2}(-k_j/2-u_0)\,\omega\wedge\omega_u^{n-2}\wedge\omega_g+2\,\int\limits_{u_0<-k_j/2}(-u)\,dd^cu_0\wedge\omega_u^{n-2}\wedge\omega_g$$$$\leq 2\,\int\limits_{u_0<-k_j/2}(-u_0)\,\omega\wedge\omega_u^{n-2}\wedge\omega_g+2\,\int\limits_{u_0<-k_j/2}(-u_0)\,\omega_{u_0}\wedge\omega_u^{n-2}\wedge\omega_g$$$$=2\,\int\limits_{u_0<-k_j/2}(-u_0)\,(\omega+\omega_{u_0})\wedge\omega_u^{n-2}\wedge\omega_g\leq 2^2\,\int\limits_{u_0<-k_j/2^2}(-u_0)\,(\omega+\omega_{u_0})^2\wedge\omega_u^{n-3}\wedge\omega_g $$$$\leq \dots\leq  2^{n-1}\,\int\limits_{u_0<-k_j/2^{n-1}}(-u_0)\,(\omega+\omega_{u_0})^{n-1}\wedge\omega_g,$$
which, by Lemma 4 and the equality $(\omega+\omega_{u_0})^{n-1}=\sum\limits_{l=0}^{n-1}{n-1\choose l}\,\omega^l\wedge\omega_{u_0}^{n-1-l}$, tends to zero as $j\to\infty$. This concludes the proof of Lemma 6.
\bigskip
We are now in a position to prove the convergence theorem.
\bigskip
\noindent {\bf Theorem 5.}(Convergence Theorem) \it Let $0\leq p<\infty$. Suppose that $0\geq g\in PSH(X,\omega)\cap L^\infty(X)$ and $u_0\in {\cal F}(X,\omega)$. If $u_j,\,u\in PSH^{-1}(X,\omega)$ are such that $u_j\to u$ in $Cap_\omega$ on $X$ and $u_j\geq u_0$, then
$(-g)^p\,\omega_{u_j}^n\to (-g)^p\, \omega_u^n$ weakly in $X$.  
\rm
\bigskip
\noindent{\it Proof.}  Given $k\geq 1$.
Write 
$$(-g)^p\,\omega_{u_j}^n-(-g)^p\,\omega_u^n=(-g)^p\,\bigl(\omega_{u_j}^n-\omega_{\max(u_j,-k)}^n\bigr)+
(-g)^p\,\bigl(\omega_{\max(u_j,-k)}^n-\omega_{\max(u,-k)}^n\bigr)$$
$$+(-g)^p\,\bigl(\omega_{\max(u,-k)}^n-\omega_u^n\bigr):=A_{k,j}+B_{k,j}+C_k.
$$
For each fixed $k$, by Theorem 1 in [X3] we have that $B_{k,j}\to 0$ weakly in $X$ as $j\to\infty$.
Given a smooth function $\psi$ in $X$. 
Following the proof of Theorem 1 in [X3], we can assume that $\psi\,(-g)^p$ is the sum of finite terms of form 
$\pm f$, where $f$ are bounded quasi-psh functions in $X$. For such a function $f$, by Lemma 2 we get
$$\Bigl|\int_Xf\,\bigl(\omega_{u_j}^n-\omega_{\max(u_j,-k)}^n\bigr)\Bigr|=\Bigl|\int_X(u_j-\max(u_j,-k))\,dd^cf\wedge \sum\limits_{l=0}^{n-1}\omega_{u_j}^l\wedge\omega_{\max(u_j,-k)}^{n-1-l}\Bigr|$$ $$=\Bigl|\int_{u_j<-k}(u_j+k)\,dd^cf\wedge \sum\limits_{l=0}^{n-1}\omega_{u_j}^l\wedge\omega_{\max(u_j,-k)}^{n-1-l}\Bigr|\leq\int_{u_j<-k}(-u_j)\,(\omega_f+\omega)\wedge \omega_{u_j}^{n-1},$$ which by Lemma 6 tends to zero uniformly for all $j$ as $k\to\infty$. Hence, $A_{k,j}\to 0$ uniformly for all $j$ as $k\to\infty$. Similarly, we have that $C_k\to 0$ weakly as $k\to\infty$. Therefore, we have obtained that $(-g)^p\,\omega_{u_j}^n\to (-g)^p\, \omega_u^n$ weakly and the proof of Theorem 5 is complete.
\bigskip
Applying Dini's theorem and quasicontinuity of quasi-psh functions, we get the following consequence.
\bigskip
\noindent {\bf Corollary 5. } \it Let $0\leq p<\infty$ and $0\geq g\in PSH(X,\omega)\cap L^\infty(X)$. If $u_j,\,u\in {\cal F}(X,\omega)$ are such that $u_j\searrow u$ or $u_j\nearrow u$ in $X$, then
$(-g)^p\,\omega_{u_j}^n\to (-g)^p\,\omega_u^n$ weakly in $X$.  
\rm
\bigskip
\noindent {\bf Corollary 6. } \it Let $u,\,v\in {\cal F}(X,\omega)$. Then
$$\chi_{\{u>v\}}\,\omega_{\max(u,v)}^n=\chi_{\{u>v\}}\,\omega_u^n.$$
\rm
\bigskip
\noindent{\it Proof.} This proof is similar to the proof of Theorem 4.1 in [KH]. Given a constant $k\geq 0$. Write $u_j=\max(u,-j)$. By Proposition 4.2 in [BT1] we have that $\max(u_j+k,0)\,\omega_{\max(u_j,-k)}^n=\max(u_j+k,0)\,\omega_{u_j}^n$ for all $j$. Using $\max(u_j+k,0)\geq \max(u+k,0)\geq 0$, we get $\max(u+k,0)\,\omega_{\max(u_j,-k)}^n=\max(u+k,0)\,\omega_{u_j}^n$. Letting $j\to\infty$ and applying Theorem 5, we get $\max(u+k,0)\,\omega_{\max(u,-k)}^n=\max(u+k,0)\,\omega_u^n$. Hence we have obtained that $\chi_{\{u>-k\}}\,\omega_{\max(u,-k)}^n=\chi_{\{u>-k\}}\,\omega_u^n$ holds for any $u\in {\cal F}(X,\omega)$ and $k\geq 0$.
Therefore, $\omega_{\max(u,v)}^n=\omega_{\max(u,v,-k)}^n$ and $\omega_u^n=\omega_{\max(u,-k)}^n$  on each set $\{u>-k>v\}$ with a rational number $k\geq 0$. But $\omega_{\max(u,v,-k)}^n=\omega_{\max(u,-k)}^n$ on the open set $\{-k>v\}$ and hence $\chi_{\{u>-k>v\}}\,\omega_{\max(u,v)}^n=\chi_{\{u>-k>v\}}\,\omega_u^n$, which implies the required equality. The proof of Corollary 6 is complete.
\bigskip
\noindent {\bf Corollary 7. } \it Let $u,\,v\in {\cal F}(X,\omega)$. Then
$$\omega_{\max(u,v)}^n\geq \chi_{\{u\geq v\ {\rm and\ } u\not=-\infty \}}\,\omega_u^n+\chi_{\{u<v\}}\,\omega_v^n.$$
\rm
\bigskip
\noindent{\it Proof.} Given $\varepsilon>0$, by Corollary 6 we have $\omega_{\max(u,v-\varepsilon)}^n\geq \chi_{\{u>v-\varepsilon\}}\,\omega_u^n+\chi_{\{u<v-\varepsilon\}}\,\omega_v^n\geq \chi_{\{u\geq v\ {\rm and\ } u\not=-\infty \}}\,\omega_u^n+\chi_{\{u<v-\varepsilon\}}\,\omega_v^n.$ Letting $\varepsilon\searrow 0$ and using Theorem 5, we obtain 
the required inequality which concludes the proof.
\bigskip
\noindent {\it Remark.} Corollary 7 is a generalization of the well known Demailly inequality, see [D2].
\bigskip
\noindent 
{\bf Corollary 8. } \it Let $u,\,v\in {\cal F}(X,\omega)$. Then
$$\int\limits_{u<v}\omega_v^n\leq \int\limits_{u<v}\omega_u^n+\int\limits_{u=v=-\infty}\omega_u^n.$$
\rm
\bigskip
\noindent{\it Proof.} By Corollary 6 we have $\int_{u<v}\omega_v^n=\int_{u<v}\omega_{\max(u,v)}^n=\int_X\omega^n-\int_{u\geq v}\omega_{\max(u,v)}^n\leq \int_X\omega^n-\int_{u>v}\omega_{\max(u,v)}^n=\int_X\omega^n-\int_{u>v}\omega_u^n=\int_{u\leq v}\omega_u^n$. Using $\delta\,v$ instead of $v$ and letting $\delta\nearrow 1$, we get the required inequality and the proof is complete.
\bigskip

\vfill\eject

\bigskip \bigskip
\bigskip \centerline{\bf References } \bigskip
\bigskip 
\noindent [B] $\,$E.Bedford, {\it Survey of pluripotential theory}. Several complex variables(Stockholm 

1987/1988) Mathematical Notes 38, Princeton University Press, 1993, pp.48-97.

\noindent [BL1] $\,$Z.Blocki, {\it Uniqueness and stability of the complex Monge-Amp\`ere operator on 

compact K\"ahler manifolds}. Indina Univ. Math. J.   {\bf 52} (2003), no 6, 1697-1701. 

\noindent [BL2] $\,$Z.Blocki, {\it The domain of definition of the complex Monge-Amp\`ere operator}. Amer. 

J. Math.  {\bf 128} (2006), 519-530. 

\noindent [BT1] $\,$E.Bedford and B.A.Taylor, {\it Fine topology, \v Silov boundary and
$(dd^c)^n$}. J. Funct. 

Anal., {\bf 72} (1987), 225-251.

\noindent [BT2] $\,$E.Bedford and B.A.Taylor, {\it A new capacity for plurisubharmonic
functions}. Acta 

Math., {\bf 149} (1982), 1-40.

\noindent [C1] $\,$U.Cegrell, {\it
Pluricomplex energy}. Acta Math. {\bf 180:2} (1998), 187-217.

\noindent [C2] $\,$U.Cegrell, {\it The general definition of the complex Monge-Amp\`ere operator }. Ann. Inst. 

Fourier {\bf 54} (2004), 159-179.

\noindent [D1] $\,$J.-P.Demailly, {\it Regularization of closed positive currents and intersection theory}. J. 

Alg. Geom., {\bf 1} (1992), 361-409.

\noindent [D2] $\,$J.-P.Demailly, {\it Potential theory in several complex variables}. \'Ecole d'\'et\'e d'Analyse 

complexe du CIMPA, Nice, 1989.

\noindent [GZ] $\,$V.Guedj and A.Zeriahi, {\it The weighted Monge-Amp\`ere energy of quasiplurisubhar-

monic functions.} , available at www.arxiv.org: math.CV/061230.

\noindent [KH] $\,$N.V.Khue and P.H.Hiep, {\it Some properties of the complex Monge-Amp\`ere operator in 

Cegrell's classes and applications.} Preprint, available at www.arxiv.org:  0704.0359.

\noindent [KI]  $\,$C.O.Kiselman, {\it Sur la definition de l'op\'erateur de 
Monge-Amp\`ere complexe.} Analyse 

Complexe: Proceedings, Toulouse (1983), 139-150. LNM 1094. Springer-Verleg.

\noindent [KO1] $\,$ S.Kolodziej, {\it The complex Monge-Amp\`ere equation and pluripotential theory } 

Memoirs of the Amer. Math. Soc. Vol.178, No 840. 2005.

\noindent [KO2] $\,$ S.Kolodziej, {\it The set of measures given by bounded solutions of the complex 
Monge-

Amp\`ere equation on compact K\"ahler manifolds.} J. London Math. Soc., (2) {\bf 72} (2005), 

225-238.

\noindent [X1] $\,$Y.Xing, {\it Continuity of the complex Monge-Amp\`ere operator.} 
 Proc. of Amer. Math.  

Soc., {\bf 124} (1996), 457-467.

\noindent [X2] $\,$Y.Xing, {\it Convergence in capacity.} Ume\aa\ University, Research Reports, No 11, 2006.

\noindent [X3] $\,$Y.Xing, {\it Continuity of the complex Monge-Amp\`ere operator on compact K\"ahler mani-

folds.} Ume\aa\ University, Research Reports, No 3, 2007, also available at www.arxiv.org: 

math/0703755.

\end